\documentclass{scrartcl}
\usepackage[T1]{fontenc}
\usepackage[latin1]{inputenc}
\usepackage{a4,amssymb,amsmath}
\usepackage[mathscr]{eucal}
\usepackage{enumerate}
\usepackage{color}

\newtheorem{thm}{Theorem}[section]

\newtheorem{lem}[thm]{Lemma}

\newtheorem{defn}[thm]{Definition}

\newtheorem{rem}[thm]{Remark}

\newtheorem{ack}[thm]{Acknowledgement}

\newcommand{\norm}[1]{\left\Vert#1\right\Vert}
\newcommand{\abs}[1]{\left\vert#1\right\vert}

\newcommand{\Cp}{\mathbb C}
\newcommand{\eps}{\varepsilon}

\newcommand{\C}{\mathcal{C}}
\newcommand{\sm}{\smallskip \\}

\newcommand{\db}{\overline\partial}

\newcommand{\be}{\begin{equation}}
\newcommand{\ee}{\end{equation}}
\newcommand{\bes}{\begin{equation*}}
\newcommand{\ees}{\end{equation*}}

\title{Exposing points on the boundary of a strictly pseudoconvex or a locally convexifiable domain of finite 1-type}%

\author{K. Diederich, J. E. Forn\ae ss and E. F. Wold\footnote{Supported by the NFR grant 209751/F20}}

\date{\today}

\begin{document}
\maketitle

\noindent \begin{abstract} \textsc{Abstract}:
We show that for any bounded domain $\Omega\subset\Cp ^n$ of 1-type $2k $ which is locally convexifiable at $p\in b\Omega$, having a Stein neighborhood basis,  there is a biholomorphic map $f:\bar{\Omega}\rightarrow \Cp ^n $ such that $f(p)$ is a global extreme point of type $2k$ for $f{(\overline\Omega)}$.
\end{abstract}

\section{Introduction}
In this paper we consider bounded locally convexifiable domains
$\Omega$ in $\mathbb{C}^n$ of finite 1-type whose closures $\bar\Omega$ admit a Stein neighborhood basis.
Here the term "locally convexifiable near $p\in b\Omega$" means that
there are a neighborhood $V$ of $p$ and a one-to-one holomorphic map $\Phi: V\to\mathbb{C}^n$
such that $\Phi(\Omega\cap V)$ is convex.  For the notion of finite type we refer to \cite{Dangelo}.
Strongly pseudoconvex domais are examples of such domains.
We will first prove the following:

\begin{thm}\label{T:main} Let $\Omega\subset\Cp^n$ be a bounded domain which is locally convexifiable and has finite type $2k$
near a point $p\in b\Omega$.  Assume further that $b\Omega$
is $\C^{\infty}$-smooth near $p$, and that $\overline\Omega$ has a Stein neighborhood basis. Then there exists a holomorphic
embedding $f:\overline{\Omega}\rightarrow \overline{\mathbb{B}}_k ^n$, where $\mathbb{B}_k ^n=\{z\in\Cp ^n:\abs{z_n}^2+\|z'\|^{2k}<1\}$, such that
\begin{enumerate}

\item $f(p)=e_n=(0,\dots,0,1)$, and
\item $\{z\in \overline{\Omega}: f(z)\in b\mathbb{B}^n_k\} = \{p\}$.
\end{enumerate}
If $k=1$, \emph{i.e.}, if $b\Omega$ is strongly pseudoconvex near $p$, it is enough to assume that $b\Omega$ is
$\mathcal C^2$-smooth near $p$.
\end{thm}

\begin{defn}
Let $\Omega\subset\mathbb C^n$ be a domain and let $p\in b\Omega$
be a point.  We say that $p$ is a globally exposed $2k$-convex
point if there exists an affine linear map $f$ as in the previous theorem.
\end{defn}

One of our motivations for proving this theorem is the special case of strictly pseudoconvex domains. In this case the theorem answers a question posed by Fusheng Deng (private communication), and it is a step to study squeezing functions on bounded strongly pseudoconvex domains (see \cite{Deng2012}). \\
Another motivation consists in the construction of
a $\C^{\infty}$ family of local holomorphic support functions $\hat{S}(z,\zeta)\in \C^{\infty}(\Cp ^n \times \partial D)$ for locally convexifiable domains of finite type $2k$ with Stein neighborhood basis as explained in \cite{Diederich1999}. It has been asked several times whether these support functions can always be chosen such that they are globally supporting for the given domains. \\
However, it has to be asked in which way precisely this question should be answered.  As far as we can see, there are at least the following different possibilities, each of them leading to quite different answers:
\begin{enumerate}
\item Our original support surfaces are defined only locally. The danger might be, for instance, that after a while they fall back into the inside of the domain or, at least, become tangent at certain points, that are further away. However, this danger can be avoided by applying a simple standard $\db$-argument to the defining functions of the support functions. Then we get new support functions which are well-defined in a possibly narrow Stein neighborhood of $\overline{\Omega}$.
\item Asking for more might mean that we really want \emph{globally} defined support surfaces, \emph{i.e.}, support surfaces which are closed smooth complex hypersurfaces  in  $\Cp^n$, touching $b\Omega$ only from the outside at one distinguished boundary point.  It is clear that this requires a much stronger hypothesis on the domain.   Namely, we will assume that the given domain has a Runge neighborhood basis and is locally convexifiable of finite type near $0$. It is one of the main results of this article (Theorem \ref{T:3})  that such closed global support surfaces then always do exist.  Under suitable regularity assumptions on $b\Omega$ (namely $b\Omega$ has to be $\C^{\infty}$-smooth) smooth $\C^\infty$-families of such supporting hypersurfaces do indeed exist (Theorem \ref{T:parameters}).
\end{enumerate}
In this  part of the work we will prove the following statement:

\begin{thm}\label{T:3}
Assume in addition to the hypotheses in Theorem \ref{T:main} that $\overline{\Omega}$ has a Runge and Stein neighborhood.
Then the map $f$ can be chosen as a global automorphism of $\Cp^n$.  A special case of this are convex domains of finite 1-type.
\end{thm}

Finally, in the case of bounded and smooth convex domains, we prove a version
of Theorem \ref{T:main} with parameters:

\begin{thm}\label{T:parameters}
Let $\Omega\subset\mathbb C^n$ be a smooth and bounded convex
domain of finite type $2k$.  There exists a smooth
family $\psi_{\zeta}\in Aut_{hol}\Cp ^n$, $\zeta\in b\Omega$,
such that $\psi_\zeta(\zeta)$ is a globally exposed $2k$-convex boundary point
for the domain $\psi_\zeta(\Omega)$.\end{thm}

\medskip

The structure of the article is as follows: In Section 2 we recall some local properties of convexifiable domains due
to the two first authors.  In Section 3 we prove Theorem \ref{T:main}.
In Section 4 we prove Anders\'{e}n-Lempert theorems with parameters needed to prove Theorem \ref{T:parameters},
which we will do in Section 5.  Finally, in Section 6, we give a brief sketch of how to prove Theorem \ref{T:3} based on the arguments in Sections 3 and 5.

\section{Local properties of convexifiable domains}
Let $\Omega$ be a bounded $\mathcal C^\infty$-smooth domain in $\mathbb C^n$.
In this section we recall the main facts about supporting hypersurfaces constructed in \cite{Diederich1999}
For this we suppose that there is an open set $V\subset\mathbb C^n$ such that  $b\Omega\cap V$ is convex.
Near any point $\zeta_0\in b\Omega\cap V$ there is an open neighborhood $V_{\zeta_0}$ of $\zeta_0$, and a choice of a $\C ^{\infty}$-family of coordinate changes $\{l_\zeta (z): \zeta \in b\Omega\cap V_{\zeta_0}\}$ composed of a translation and a unitary transformation, such that, for each $\zeta \in b\Omega\cap V_{\zeta_0}$, $l_\zeta (\zeta)=0$
and the unit outward normal vector $n_\zeta$ at $\zeta$ is turned by $l_\zeta$ into the unit vector $(1,0,\ldots,0)$. In particular, $T^{\mathbb C}_\zeta b\Omega$ becomes in the new coordinates $\tilde z=l_\zeta (z)$ associated to $\zeta$ just $\{\tilde z_1=0\}$.
The following is proved in \cite{Diederich1999}:

\begin{thm}
In the situation just described, assume that $b\Omega\cap V$ is of finite 1-type $2k$, and let $\tilde V\subset\subset V$.
Then there exists a function $\widehat S(\zeta,z)\in\mathcal C^\infty((b\Omega\cap\tilde V)\times\mathbb C^n)$,
and constants $r,c>0$,
such that the following holds:  for any choice of coordinate changes $l_\zeta$ as above, the
function $S(\zeta,z):=\widehat S (\zeta,l_\zeta^{-1}(z))$ is equal to
\begin{equation}
S_\zeta(z)=3z_1 + Kz_1^2 + g_\zeta(z'), \mbox{ where } (z_1,z') \mbox{ are coordinates on } \mathbb C^n,
\end{equation}
and satisfies the estimate
\begin{equation}\label{estimate}
Re S_\zeta(z)\leq -c(|z_1|^2+\|z'\|^{2k}),
\end{equation}
for all $z\in B_r\cap l_\zeta(\Omega)$.
\end{thm}

Note that if the domain $\Omega$ is convex, we get a $\mathcal C^\infty$-smooth function $\widehat S(\zeta,z)$
on $b\Omega\times\mathbb C^n$.\sm


\section{The proof of Theorem \ref{T:main}}

The proof of Theorem 1.1 is reduced to the two Lemmas in this section,
Lemma \ref{L:1} and Lemma \ref{L:3}.  \\

We let ${e_1,\ldots,e_n}$ be the standard basis for the complex vector space $\Cp^n$ and put  $f_j:=i\cdot e_j$ so that ${e_1,f_1,\ldots,e_n,f_n}$ is a real basis. We denote the coordinates on $\Cp ^n$ by $z_j =
x _{2j-1}+ix _{2j}$, we let $\Cp_n$ denote the complex line $\Cp _n =\{z\in \Cp ^n : z_1 =\ldots =z _{n-1}=0\}$ and we let $\pi _n$ be the orthogonal projection to $\Cp _n$.\\
Our proof uses a technique from \cite{Forstneric2009.1} invented for exposing points on a bordered Riemann surface in order to produce a proper holomorphic embedding (see also \cite{Forstneric2011} Sections 8.8 and 8.9).
We suppose that $\Omega$ is convexifiable near some point $p$ on its boundary. Then we get the following situation:
\begin{lem}\label{L:1}
For any $p\in b\Omega$ there exists $\Phi \in {\rm Aut} _{hol}\Cp ^n$ such that the following hold
\begin{enumerate}
\item $\Phi (p)=0$ and $T_0 (b\Phi(\Omega))=\{x_{2n}=0\}$\label{L:1-1}
\item The outward normal to $b\Phi(\Omega)$ at the origin is $f_n$,\label{L:1-2}
\item Near the origin we have that $b\Phi {(\Omega)}$ is $k$-convex at the origin in the following sense: The domain $\overline{\Phi (\Omega)}\subset \{z\in \Cp^n : x_{2n}-f(z',x_{2n-1}) \leq 0\}$ with $f(z',x_{2n-1})\geq c(\|z'\|^{2k}+x^{2}_{2n-1})$, $c>0$ and\label{L:1-3}
\item $\overline{\Phi (\Omega)}\cap\{z\in \Cp^n: z_1 =\ldots= z_{n-1} =x_{2n-1}=0,x_{2n}\geq 0\}=\{0\}$\label{L:1-4}

\end{enumerate}
\end{lem}
\begin{defn}
When condition (3) is satisfied near the origin we will refer to the origin as a strictly $2k$-convex boundary point.
\end{defn}

\underline{Proof:}
It follows by Corollary 2.4 in \cite{Diederich1999} that there exists an open neighborhood $U_p$ of $p$ and an injective holomorphic map $\psi :U_p \rightarrow \Cp ^n$ such that $\psi(\Omega \cap U_p)$ satisfies (1)-(3). Choosing an appropriate neighborhood $V_p\subset U_p$ of $p$, it follows  that $\psi$ is approximable by automorphisms $\phi$ of $\Cp ^n$ uniformly on $V_p$ (see Section 4) and that $\psi(\Omega)$ is strictly k-convex near $\Phi (p)$ if $\phi$ is close enough to $\psi$.\\
We proceed to achieve (4). Let
\be
\Gamma :=\{z\in\Cp ^n :z_1 =\ldots =z_{n-1} = x_{2n-1}=0, x_{2n}\geq 0\},
\ee
and let
\be
\Gamma _0:=\{z\in \Cp^n: z_1=\ldots =z_{n-1}=x_{2n-1}=0, 0\leq x_{2n}\leq 1\}.
\ee
Choose an $R>0$ such that $\overline{\Omega}\subset \mathbb{B}^n_R$. By \cite{Forstneric1997b} there exists $\psi _1 \in
{\rm Aut}_{hol}(\Cp ^n)$ such that $\psi_1 (z)=z+O(\norm{z}^{2k+1})$ as $z\rightarrow 0$, such that $\psi _1 (\Gamma _0)\cap\overline{\Omega}= \{0\} $, and such that $\psi _1(q)\in \Cp ^n\backslash \overline{\mathbb{B}^n _R}$,where $q$ denotes the endpoint of $\Gamma _0$ other than the origin. Consider the set $\psi ^{-1} _1 (\overline{\mathbb{B}^n _R}) \cap \Gamma \subset \Gamma \backslash \{q\}$ . Since
$\psi ^{-1}(\overline{\mathbb{B}^n _R})$ is polynomially convex we have that $(\Cp ^{n-1} \times \{i\})\backslash \psi ^{-1} _1(\overline{B ^n _R})$ is connected, and so using Weierstrass approximation theorem, we may construct a holomorphic shear map $\psi _2 (z)=(z_1+f_1(z),\ldots,  z_{n-1}+f_{n-1}(z_n), z_n)$ such that
$\psi_2$ is close to the identity on $\Gamma_0$, tangent to the identity to order 2k+1 at the origin, and therefore not destroying strict k-convexity at 0, and such that $\psi_2 (\overline{(\Gamma \backslash \Gamma_0)}\cap \psi^{-1} _1 (\overline{B}^n _R) = \emptyset $. So $(\psi_1\circ \psi_2)(\Gamma)\cap\overline{\Omega}=\emptyset $, and we set $\Phi =(\psi _1 \circ \psi _2)^{-1}$. \sm

\begin{lem}\label{L:3}
Let  $W \subset \Cp ^n$ be a bounded domain with $0\in bW$ and assume that the following hold
\begin{enumerate}[\upshape (i)]
\item $\overline{W}$ has a Stein neighborhood basis,
\item $W$ is strictly k-convex near the origin,
\item $\overline{W}\cap\Gamma _0 = \{0\}$ with $\Gamma _0$ defined as in the proof of Lemma 3.1.
\end{enumerate}
Then for any open set $\tilde{V}$ containing $\Gamma _0$ and any small enough open set $V\subset \tilde{V}$ containing the origin, there exist a sequence of holomorphic embeddings $f_j:\overline{W} \rightarrow \Cp ^n$ such that the following holds
\begin{enumerate}[\textbf (1)]
\item $f_j\rightarrow id$ uniformly on $\overline{W}\setminus V$ as $j\rightarrow \infty$\label{L:3-1}
\item $f_j (V)\subset \tilde{V}$ for all $j$,\label{L:3-2}
\item $f_j (0)=f_n$ for all $j$,\label{L:3-3}
\item ${\rm Im}(\pi_n (f_j(z))) < 1$ for all $z\in (\overline{W}\cap V)\setminus \{0\}$ and\label{L:3-4}
\item $f_j(W)$ is strictly k-convex at $f_j (0)$.\label{L:3-5}
\end{enumerate}
\end{lem}
\underline{Proof:} We let $\tilde{\Omega}\subset \Cp _n$ denote the domain $\tilde{\Omega}:= \{z_n\in\Delta_\eps : x_{2n} <f(0,\ldots,0,x_{2n-1}\})$. For some small $\delta >0$ we define the following sets: $A:=\{z\in \overline{W}\cap\overline{\mathbb{B}^n _{\eps}}\,: x_{2n}\geq -\delta $ and $B:=\{z\in\overline{W}\cap\overline{\mathbb{B}^n _{\eps}}\,:x_{2n}\leq -\frac{\delta}{2}\}\cup\overline{W}\setminus \mathbb{B}^n _{\eps}$. Then  $\overline{A\setminus B}\cap\overline{B\setminus A}=\emptyset$ (if $\delta$ is small) and by Theorem 4.1 in \cite{Forstneric2003}, for any open set $\tilde{C}$ containing the set $C:=A\cap B$ there exist open sets $A',B',C'$ with $A\subset A'$, $B\subset B'$ and $C\subset C'\subset{A'\cap B'}\subset \tilde{C}$, such that if $\gamma :\tilde{C}\rightarrow \Cp ^n$ is injective holomorphic, and sufficiently close to the identity, then there exist holomorphic injections $\alpha : A'\rightarrow \Cp ^n$, $\beta : B'\rightarrow \Cp ^n$, uniformly close to the identity on their respective domains (depending on $\gamma$), and such that
\be
\gamma = \beta\circ\alpha ^{-1} \text{ on $C'$.}
\ee
(This can also be found in Theorem 8.7.2, page 359 in \cite{Forstneric2011}.)
Choose a simply connected smooth domain $U\subset \Cp _n$ with $\pi _n (A)\subset U$ and such that near the origin $U=\{z\in\Cp _n :x_{2n} < 0\}$. For $j\in\mathbb{N}$ let $l_j$ denote the line segment $l_j=\{z_n\in\Cp _n : x_{2n-1} =0,0\leq x_{2n}\leq 1/j\}$. For each $j$ it follows from Mergelyan's Theorem that we may choose injective holomorphic maps $\sigma _j :\overline{U}\cup l_j\rightarrow \Cp _n$ such that $\sigma_j$ approximately stretches $l_j$ to cover $\Gamma _0$ such that $\sigma _j(z)=(1-1/j)i+z+O(\abs{z-i/j})^{2k+1}$ and such that $\sigma _j\rightarrow id$ on $\overline{U}$ as $j\rightarrow\infty$. For each $j$ let $U_j$ be a domain obtained from U by adding a strip around $l_j$ of width less than $1/j$ which is then smoothened and made strictly convex at the end point $l_j$. $U_j$ should lie inside where $\sigma _j$ is injective holomorphic, and be chosen such that $\sigma _j (U_j)$ is strictly convex near the end point of $\sigma _j (l_j)=f_n$ and such that ${\rm Im}(\sigma _j(z_n))<1$ for all $z\in \overline{U_j}\backslash \frac{1}{j} f_n$. Let $\psi_j$ be a holomorphic diffeomorphism  from $\overline{U}$ to $\overline{U}_j$ such that $\psi_j(0)=\frac{i}{j},$  $\psi _j\rightarrow id$ uniformly on $\overline{U}$. (See Goluzin, \cite{Goluzin1969}, Theorem 2, p. 59.) Let $\phi _j =\sigma _j\circ \psi _j$ and let $\gamma _j$ be an extension of $\phi _j$ to A. Then ${\rm Im}\left(\Pi_n(\gamma _j(z))\right)<1$ for all $z\in A\setminus \{0\}$. It is not hard to see that $\gamma _j(A)$ is strictly k-convex near $f_n$ and $\gamma _j\rightarrow id$ on a neighborhood of $C$. We get splittings
\be
\gamma _j\circ \alpha _j = \beta _j
\ee
as explained above. If $j$ is large enough, we get that (7) defines an injective holomorphic map $f_j$ on $\overline{\Omega}$, and if $\alpha _j$ is close enough to the identity, since $\alpha _j$ can be assumed to vanishes to order $2k+1$ at the origin, we get that ${\rm Im}(f_j (z))<1$ for all $z\in A\setminus\{0\}$ and such that $f _j(A)$ is strictly k-convex at $f _j(0)$.


\section{Anders\'{e}n-Lempert with parameters in a smooth manifold, and approximation with
jet interpolation.}

A parameter version of the Andersen-Lempert theorem \cite{Andersen1992} for holomorphic parameters was proved by Kutzschebauch\cite{Kutzschebauch2005}.
Jet interpolation results without parameters have been proved by Forstneri\v{c} \cite{ForstnericJet} and Weickert \cite{Weickert}
(see also sections 4.9 and 4.15 in \cite{Forstneric2011}).
For a smooth manifold $M$ we let $(\zeta,z)$ denote the coordinates on $M\times\mathbb C^n$.
For any $\zeta\in M$ we denote by $\mathbb C^n_\zeta$ the slice $\{\zeta\}\times\mathbb C^n$,
and for any subset $\Sigma\subset M\times\mathbb C^n$ we let $\Sigma_\zeta$ denote
the slice $\Sigma_\zeta:=\mathbb C^n_\zeta\cap\Sigma$.

\begin{thm}\label{alparameters}
Let $M$ be a compact smooth manifold and let $\Omega\subset M\times\mathbb C^n$ be a domain, $n\geq 2$.
Let $K\subset\Omega$ be a compact set, and let $\phi:[0,1]\times\Omega\rightarrow M\times\mathbb C^n$ be a $\mathcal C^2$-smooth
map such that, writing $\phi(t,\zeta,z)=\phi_t(\zeta,z),$ the following hold
\begin{itemize}
\item[(1)] $\phi_t(\zeta,z)=(\zeta,\varphi_t(\zeta,z))=(\zeta,\phi_{t,\zeta}(z))$,
\item[(2)] $\phi_{t,\zeta}:\Omega_\zeta\rightarrow\mathbb C^n_\zeta$ is injective holomorphic, and
\item[(3)] $K_{t,\zeta}:=\phi_{t,\zeta}(K_\zeta)$ varies continuously with $(t,\zeta)$ and is polynomially convex
\end{itemize}
for all $t\in [0,1], \zeta\in M$.  \sm
Then $\phi_1$ is uniformly approximable on $K$ by
a smooth family $\psi(\zeta,z)$ with $\psi_\zeta\in Aut_{hol}\mathbb C^n_\zeta$
if (and only if) $\phi_0$ is approximable by such a family.   Moreover, if (1)--(3) hold and if $a(\zeta)\in K^{\circ}_\zeta$
is a smoothly parametrized family of points, and if $d\in\mathbb  N$,
we may additionally achieve that
\begin{itemize}
\item[(4)] $\phi_{1,\zeta}(z)-\psi_\zeta(z)=O(\|z-a(\zeta)\|^d)$, as $z\rightarrow a(\zeta)$.
\end{itemize}
\end{thm}
\underline{Proof:}
We give a sketch of the proof of the first claim; the point is just to verify that the non-parametric
proof goes through without change with parameters.\\
The assumption that $\phi_0$ is approximable allows us to assume $\phi_0=id$.
Define first a parametrized vector field
\begin{equation}
X_{t,\zeta}(\phi_{t,\zeta}(z)):=\frac{d}{dt}\phi_{t,\zeta}(z).
\end{equation}
Then $X_{t,\zeta}$ is an inhomogeneous vector field, holomorphic in $z$, whose
flow is $\phi_{t,\zeta}(z)$.  For each $t$ let $\varphi^s_{t,\zeta}$ denote
the time-$s$ flow of the homogenous vector field $X_{t,\zeta}$ where $t$ is fixed.  It is well known that there is a partitioning
$[j/n,(j+1)/n], j=0,...,n-1$ of $[0,1]$, such that the composition
\begin{equation}
\varphi^{1/n}_{(n-1/n,\zeta)}\circ\cdot\cdot\cdot\circ\varphi^{1/n}_{0,\zeta}
\end{equation}
approximates $\phi_{\zeta,1}$ on $K$.  So the problem is reduced to approximating
the flow $\varphi^1_\zeta$ of a homogenous vector field $X_\zeta$ on a family $K_\zeta$.

\medskip

Next, by assumption (3) and approximation, we may assume that $X_\zeta$ is a polynomial vector field
\begin{equation}\label{sum}
X_\zeta(z)=\sum_{j=1}^N g_j(\zeta)X_j(z),
\end{equation}
with
coefficients $g_j$ in $\mathcal E(M)$; this can be obtained by gluing a fiberwise Runge-approximation
using a partition of unity on $M$.  Now the main point of Anders\'{e}n-Lempert Theory in $\mathbb C^n$
is that \emph{any $m$-homogenous polynomial vector field $V_m$ is a sum of completely integrable vector fields} (see e.g.
\cite{Forstneric2011}, Lemma 4.9.5):
\begin{equation}\label{shearfields}
V_m(z)= \sum_{i=1}^r c_i\lambda_i(z)^m\cdot v_i + d_i\lambda_i(z)^{m-1}\langle z,v_i\rangle\cdot v_i,
\end{equation}
with $c_i,d_i\in\mathbb C, v_i\in\mathbb C^n$ and $\lambda_i\in(\mathbb C^n)^*$ with $\lambda_i(v_i)=0$. The flows of these two types of vector fields are
\begin{equation}\label{shearflows}
z\overset{f_{t,j}}{\mapsto} z + t\cdot c_i\lambda_i(z)^m\cdot v_i \mbox{ and } z \overset{g_{t,j}}{\mapsto} z+(e^{td_i\lambda_i(z)^m}-1)\langle z,v_i\rangle\cdot v_i.
\end{equation}

Applying this to each of the vector fields $X_j(z)$ in (\ref{sum}) we get that
\begin{equation}
X_\zeta(z)=\sum_{j=1}^{\tilde N} \tilde g_j(\zeta)\cdot \tilde X_{j}(z),
\end{equation}
where each $\tilde X_{j}$ is completely integrable with flow $\psi_j^s$, and so
$X_\zeta$ is a sum of completely integrable fields
with flows $\psi_{\zeta,j}^s=\psi_j^{g(\zeta)\cdot s}$.  Finally the sequence
\begin{equation}
(\psi^{1/n}_{\zeta,N}\circ\cdot\cdot\cdot\circ\psi^{1/n}_{\zeta,1})^n
\end{equation}
converges uniformly to $\varphi^1_\zeta$ as $n\rightarrow\infty$.

\medskip

Finally we consider (4).   We will correct the initial approximation at $a(\zeta)$ and
by translation we may assume that $a(\zeta)=0$ for all $\zeta$, and that
both $\phi$ and $\psi$ fix the origin.   Define
$J_{d-1}(\zeta)$ to be the $(d-1)$-jet of $\psi_\zeta^{-1}\circ\phi_{1,\zeta}$.
It is easy to see that we may assume that $J_{d-1}(z)=id + h.o.t$, and
by the Cauchy estimates we may assume that $J_{d-1}(\zeta)$ is
arbitrarily close to the identity map.  We will correct $\psi_\zeta$
inductively, and our induction assumption is that $J_{d-1}(\zeta)=O(\|z\|^{m}), 2\leq m\leq d-1$.

\medskip

Using (\ref{shearfields}) we fix an expansion
\begin{equation}\label{expansion}
z^\alpha\cdot e_j = s_{\alpha,j}(z):=\sum_{i=1}^r c^{\alpha,j}_i\lambda^{\alpha,j}_i(z)^m\cdot v^{\alpha,j}_i + d^{\alpha,j}_i\lambda^{\alpha,j}_i(z)^{m-1}\langle z,v^{\alpha,j}_i\rangle\cdot v_i
\end{equation}
for each multi-index $|\alpha|=m$ and $j=1,...,n$.   Now expand the $m$-homogenous part $J_{d-1,m}$ of $J_{d-1}$ using (\ref{expansion})
\begin{equation}
J_{d-1,m}(\zeta)=\underset{|\alpha|=m, 1\leq j\leq n}{\sum} h_{\alpha,j}(\zeta)\cdot s_{\alpha,j}(z).
\end{equation}
It is easy to see that the composition $\Phi_m$ of all automorphisms

\begin{equation}
z\mapsto z + h_{\alpha,j}(\zeta)\cdot c^{\alpha,j}_i\lambda^{\alpha,j}_i(z)^m\cdot v^{\alpha,j}_i
\end{equation}
and
\begin{equation}
z  \mapsto z+(e^{h_{\alpha,j}(\zeta)d^{\alpha,j}_i\lambda^{\alpha,j}_i(z)^m}-1)\langle z,v^{\alpha,j}_i\rangle\cdot v^{\alpha,j}_i.
\end{equation}
matches $J_{d-1,m}$ to order $m$, and we may assume that $\Phi_m$ is as close to the identity
as we like on a compact set since all the $h_{\alpha,j}$'s can be assumed to be as small as we like.
It follows that the map $\psi_\zeta\circ\Phi_{m}$ is a small perturbation of $\psi_\zeta$
which matches $\phi_{1,\zeta}$ to order $m$.   The induction step is complete.

\begin{rem}
For a more detailed explanation of jet-completion (without parameters) the reader can consult \cite{Forstneric2011} page 154--158.
\end{rem}


\section{The construction with parameters: Proof of Theorem \ref{T:parameters}}

\begin{thm}
Let $\Omega\subset\mathbb C^n$ be a smooth and bounded convex
domain of finite type $2k$.  There exists a smooth
parameter family $\psi_{\zeta}\in Aut_{hol}\Cp ^n$, $\zeta\in b\Omega$,
such that $\psi_\zeta(\zeta)$ is a globally exposed $2k$-convex boundary point
for the domain $\psi_\zeta(\Omega)$.
\end{thm}

\underline{Proof}:
By \cite{Diederich1999} there exist $r,c>0$ and a smooth parameter family
\be
\psi_\zeta (z) \text{ defined on } \{(\zeta,z): \norm{z-\zeta}<r\}
\ee
such that $\psi(\zeta,\cdot)$
is injective holomorphic for all $\zeta$ and the following holds for all $\zeta$ (see Section 2):
let $n_\zeta$ denote the outward pointing unit normal vector to $b\Omega$ at $\zeta$, let
$l_\zeta$ be a composition of a translation and a unitary transformation
such that $l_\zeta(\zeta)=0$ and such that $n_\zeta$ is sent to the vector $(1,0,...,0)$.
Then $\tilde\psi_\zeta(z):=\psi_{\zeta}\circ l_\zeta^{-1}$ is of the form $(S_\zeta(z),z_2,...,z_n)$,
and $S$ satisfies
\begin{equation}\label{localform}
S_\zeta(z)=3z_1+Kz_1^2 + g_\zeta(z'), \, z=(z_1,z'),
\end{equation}
(See Section 2.)
Moreover, we have that
\begin{equation}
Re(S_\zeta(z))\leq -c\cdot  \left(|z_1|^{2}+ \|z'\|^{2k}\right), z\in B_r(\zeta)\cap\overline\Omega,
\end{equation}
where the constant $c>0$ does not depend on $\zeta$.

\medskip

Our first step is to change the maps $\psi_\zeta$ conveniently on
the normals $n_\zeta$, and then approximate the changed maps
by a family of holomorphic automorphisms.
Set $\Gamma_0:=\{z\in\mathbb C^n:0\leq x_{1}\leq 1, x_2=z_2=...=z_n=0\}$,
and let $h$ denote the map $h(z)=3z_1+Kz_1^2$ near $\Gamma_0$.
By changing $h$ smoothly, then finding a smooth homotopy of maps, and
finally applying Mergelyan's Theorem with parameters, we find $\delta>0$
and a smooth map
$$
\tilde h:[0,1]\times\Gamma_0(\delta)\rightarrow\mathbb C,
$$
such that the following hold
\begin{itemize}
\item[(1)] $\tilde h_0(z_1)=3z_1$,
\item[(2)] $\tilde h_t(\cdot)$ is injective holomorphic for each $t\in [0,1]$,
\item[(3)] $\tilde h_1(z_1)\approx h(z_1)$ on $B_\delta(0)$ and $(\tilde h_1-h)(z)=O(|z|^{2k+1})$ as $z\rightarrow 0$,
\item[(4)] $\tilde h_1(z)\approx 3z_1$ on $B_\delta(1)$ and $(\tilde h_1-3)(z)=O(|z-1|^{2k+1})$ as $z\rightarrow 1$,
\item[(5)] $\tilde h_1^{-1}(3\Gamma_0)\approx \Gamma_0$.
\end{itemize}

We define a homotopy modification $\widehat\psi_{\zeta,t}(z)$ of $\psi_\zeta$ by setting
\begin{equation}
\widehat S_{\zeta,t}(z):=\tilde h_t(z_1) + t\cdot g_\zeta(z') \mbox{ on } \Gamma_0(\delta).
\end{equation}
in local coordinates.

\medskip

Let $b(\zeta)$ denote the end point of $n_\zeta$ other than $\zeta$, and note that by Stolzenberg \cite{Stolzenberg1966} we may assume, by possibly having
to decrease $\delta$, that
$$
\widehat K_{\zeta,t}:=\widehat\psi_{\zeta,t}(\overline{B_\delta(\zeta)\cup n_\zeta\cup B_\delta(b(\zeta)})
$$
is polynomially convex for all $\zeta$.   By Theorem \ref{alparameters} and its proof there
exist families $G_\zeta,H_\zeta\in Aut_{hol}\mathbb C^n$ such that the following holds
\begin{itemize}
\item[(6)] $G_\zeta\approx\widehat\psi_{\zeta,1}$ on $B_\delta(\zeta)$, and $(G_\zeta-\widehat\psi_\zeta)(z)=O(\|z-\zeta\|^{2k+1})$ as $z\rightarrow\zeta$,
\item[(7)]\label{7} $H_\zeta\approx\widehat\psi_{\zeta,1}^{-1}$ on $B_\delta(0)\cup 3n_\zeta\cup B_\delta(3b(\zeta))$,
\item[(8)]\label{8} $(H_\zeta-\widehat\psi^{-1}_{\zeta,1})(z)=O(\|z-3b(\zeta)\|^{2k+1})$ as $z\rightarrow 3b(\zeta)$, and
\item[(9)] $H_\zeta\circ G_\zeta\approx id$ on $\overline\Omega$.
\end{itemize}

\medskip

Next we construct a \emph{continuous} parameter family of exposing maps $f_\zeta$ as in Lemma \ref{L:3},
where each $f_\zeta$ wraps the boundary at $G_\zeta(\zeta)$ around the normal $3n_\zeta$.
The composition $H_\zeta\circ f_\zeta\circ G_\zeta$ will globally expose the point $\zeta$
$2k$-convexly.   We will then change $f_\zeta$ to depend \emph{smoothly} on $\zeta$,
and in a final step we will approximate the family $f_\zeta$ by a smooth family
of automorphisms.

\medskip

Choose a strictly pseudoconvex neighborhood $\Omega'$ of $\overline\Omega$ close to $\Omega$
and let $\rho$ be a smooth strictly plurisubharmonic defining function for $\Omega'$ near $b\Omega'$.
For $0<r<<1$ we let $\Omega'(r):=\{z:\rho(z)<r\}$.
For $0<\sigma<<1$ we define Cartan
pairs $\tilde A_\zeta(r):=\Omega'(r)\cap \overline B_\sigma(\zeta), \tilde B_\zeta(r):=\Omega'(r)\setminus B_{\sigma/2}(\zeta)$.
Set $A_\zeta(r):=G_\zeta(\tilde A_\zeta(r)), B_\zeta(r):=G_\zeta(\tilde B_\zeta(r))$.
\medskip

Let $\gamma_j$ be the sequence of locally exposing maps from the proof of Lemma \ref{L:3}.
Since the maps only depend on the normal coordinate, the
map $\tilde\gamma_{j,\zeta}:=l_\zeta^{-1}\circ\gamma_j\circ l_\zeta$ is a well defined family
of locally exposing maps for $\Omega$, and $\tilde\gamma_{j,\zeta}\rightarrow id$
uniformly on $C_\zeta(r):=A_\zeta(r)\cap B_\zeta(r)$ for small enough $r$ independently of $\zeta$.
To globalize these locally defined maps we use
the following parametric version of Theorem 8.7.2 in \cite{Forstneric2011}.

\begin{lem}\label{cartan}
If $r_0$ is small enough and $\mu>0$ there exist $r_1<r_0$ and $\epsilon>0$
such that the following
holds:  for any family $\gamma_\zeta:\overline C_\zeta(r_0)\rightarrow\mathbb C^n$
of holomorphic maps with $\|\gamma_\zeta-id\|_{\overline C_\zeta(r_0)}<\epsilon$, continuous in $\zeta$,
there exist injective holomorphic maps $\alpha_\zeta:A_\zeta(r_1)\rightarrow\mathbb C^n, \beta_\zeta:B_\zeta(r_1)\rightarrow\mathbb C^n$, continuous in $\zeta$, such that
\begin{equation}
\gamma_\zeta=\beta_\zeta\circ\alpha_\zeta^{-1},  \|\alpha_\zeta-id\|_{A_\zeta(r_1)}<\mu,  \|\beta_\zeta-id\|_{B_\zeta(r_1)}<\mu.
\end{equation}
Moreover, we may achieve that $(\alpha_\zeta-id)(z)=O(\|z-\zeta\|^{2k+1})$ as $z\rightarrow\zeta$.
\end{lem}
The proof of this is almost identical to that in \cite{Forstneric2011}, noting
that there exist a solution operator to the $\overline\partial$-equation which
is continuous with parameters, and one can multiply by powers of $S_\zeta$
to get exact jet interpolation.

So if $j$ is chosen large enough we get that the family $f_\zeta$ defined
as $f_\zeta:=\gamma_{\zeta,j}\circ\alpha_{j,\zeta}$ on $A_\zeta(r_1)$ and
$f_\zeta:=\gamma_{\zeta,j}\circ\beta_{\zeta,j}$ on $B_\zeta(r_1)$, is a family
of injective holomorphic maps $\tilde\gamma_\zeta:G_\zeta(\Omega'(r_1))\rightarrow\mathbb C^n$
exposing the point $\zeta$ $2k$-convexly.   By (\ref{7}) and (8) the family $H_\zeta\circ f_\zeta\circ G_\zeta$
is a continuous family of holomorphic injections on $\Omega'(r_1)$, globally exposing the point $\zeta$
for the domain $\Omega$.

\medskip

Next we approximate $f_\zeta$ by a \emph{smooth} family of exposing maps.  This is
done using a partition of unity on $b\Omega$.    Note first that although
$f_\zeta$ is only continuous in $\zeta$, the $2k$-jet at $\zeta$, $J(\zeta)$,
is smooth in $\zeta$; this is because $\alpha_{j,\zeta}$ vanishes to order $2k$
at $\zeta$.   Let $(U_j,\alpha_j), j=1,...,m$, be a partition of unity on $b\Omega$ with a point $a_j\in U_j$
for all $j$.   For each $j$ write $f_{a_j}(z)=z + g_j(z)$.  We set
$\tilde f_\zeta(z):=z + \sum_{j=1}^m\alpha_j(\zeta)g_j(z)$.  By choosing the
covering fine enough we may achieve that $\tilde f_\zeta$ is as
close to $f_\zeta$ as we like on $\overline\Omega$, and also that the $2k$-jet
of $\tilde f_\zeta$ at $\zeta$ is as close to that of $f_\zeta$ as we like.  So using
the argument in the proof of Theorem \ref{alparameters} we can correct $\tilde f_\zeta$
so that its $2k$-jet at $\zeta$ matches that of $f_\zeta$ exactly.

\medskip

Finally we need to approximate the family $\tilde f_\zeta$ by a family of automorphisms.
We may assume that $0\in\Omega, G_\zeta(0)=0$,  and that $f_\zeta(0)=0$ for all $\zeta$.
Set
\begin{equation}
\varphi_t({\zeta,z}):=G_\zeta(\frac{1}{t}\tilde f_\zeta(t\cdot G_\zeta^{-1}(z))).
\end{equation}
We may assume that $\tilde f_\zeta(G_\zeta(\overline\Omega))$
is polynomially convex for all $\zeta\in b\Omega$.  In that case
it follows that there exists some $s>1$ such that $\varphi_{t,\zeta}(G_\zeta(s\overline\Omega))$
is polynomially convex for all $t,\zeta$, and so approximation follows by Theorem \ref{alparameters}.
It is enough to show that $f_\zeta(G_\zeta(\overline\Omega))$ is
polynomially convex.

\medskip

Fix $\zeta\in b\Omega$.  By  Stolzenberg \cite{Stolzenberg1966} we have that
$G_\zeta(\overline\Omega)\cup 3n_\zeta$ is polynomially convex.   Let
$W_\zeta$ be a Runge neighborhood of $K_\zeta:=G_\zeta(\overline\Omega)\cup 3n_\zeta$,
very close to $K_\zeta$.   Consider a point $b\in b\Omega\cap\overline B_\zeta(0)$.
If $W_\zeta$ is close enough to $K_\zeta$, and if $\beta_\zeta$ is close
enough to the identity, then the locally defined function $e^{C\cdot S_b(\beta_\zeta^{-1}(z))}$
for $C>>0$ may be approximately globalized to $W_\zeta$, separating points on $\beta_\zeta(n_b)$
close to $\beta_\zeta(b)$ from $f_\zeta(G_\zeta(\overline\Omega))$ as
long as $f_\zeta$ is chosen such that $f_\zeta(\overline\Omega)\subset W_\zeta$.   It follows that
\begin{equation}
cl[\widehat{f_\zeta(G_\zeta(\overline\Omega))}\setminus f_\zeta(G_\zeta(\overline\Omega))]\cap f_\zeta(G_\zeta(\overline\Omega))\subset f_\zeta(A_\zeta(0)).
\end{equation}
Hence by Rossi's local maximum principle
\begin{equation}
\widehat{f_\zeta(G_\zeta(\overline\Omega))}=f_\zeta(G_\zeta(\overline\Omega))\cup\widehat{[f_\zeta(A_\zeta(0))\cap G_\zeta(\overline\Omega)]}.
\end{equation}
But $f_\zeta^{-1}$ is approximable by entire maps on $f_\zeta(A_\zeta(0))$, and so $f_\zeta(G_\zeta(\overline\Omega))$
is polynomially convex.

\section{Remark on the proof of Theorem \ref{T:3}}

The proof of Theorem \ref{T:3} is almost the same as that of Theorem \ref{T:main},
except that we need to make sure that the exposing maps $f_j$ are
approximable by holomorphic automorphisms.  To see why this is so, note first
that each $\gamma_j$ may be connected to the identity map by an isotopy which
is uniformly close to the identity on $C$.
The Cartan type splitting with parameters then allows us to construct each $f_j$
as the time-1 map of an isotopy $f_{j,t}$ with $f_{j,0}=\mathrm{id}$ (this argument
allows us to avoid the usual assumption in Anders\'{e}n-Lempert theory that $\Omega$
is star shaped).  This isotopy is only $\mathcal C^0$ but we
can obtain a smooth isotopy by gluing as before.  The same argument as
in the previous section tells us that we may assume that $f_{t,j}(\overline\Omega)$
is polynomially convex for all $t$ if $j$ is sufficiently large, and so we may approximate
by automorphisms.

\begin{ack}
The authors would like to thank the referee for a very carefull reading, thus
helping us to improve greatly the exposition of the paper.
\end{ack}

\bibliographystyle{amsplain}


\end{document}